\documentclass[letterpaper, 10 pt, conference]{ieeeconf}

\IEEEoverridecommandlockouts
\usepackage{cite}
\usepackage{amsmath,amssymb,amsfonts}
\usepackage{algorithmic}
\usepackage{graphicx}
\usepackage{textcomp}
\usepackage{xcolor}
\usepackage[makeroom]{cancel}
\usepackage[T1]{fontenc}

\let\labelindent\relax

\usepackage{enumerate}
\usepackage{enumitem}

\usepackage{bbm}
\usepackage{amsthm}
\usepackage{amsmath}
\usepackage{thmtools}
\usepackage{amsfonts}
\usepackage[thinc]{esdiff}
\usepackage{url}
\usepackage{hyperref}
\usepackage{cleveref}

\declaretheorem[]{lemma}
\declaretheorem[name = Theorem]{thm}

\declaretheorem[name = Assumption]{assump}

\newtheorem{dfn}{Definition}

\newcommand{\R}{\mathbb{R}}

\newcommand{\db}{\delta_{\infty}}
\newcommand{\ve}{\varepsilon}
\newcommand{\tx}{\Tilde{X}}
\newcommand{\nf}{\nabla f}
\newcommand{\nb}{\nabla B}
\newcommand{\nft}{\nabla f^{\top}}

\newcommand{\N}{\mathbb{N}}
\newcommand{\E}{\mathbb{E}}

\newcommand{\C}{\Hat{c}}
\newcommand{\Ef}{\mathbb{E}^{\mathcal{F}}_k}

\makeatletter
\newcommand{\distas}[1]{\mathbin{\overset{#1}{\kern\z@\sim}}}%
\newsavebox{\mybox}\newsavebox{\mysim}

\def\BibTeX{{\rm B\kern-.05em{\sc i\kern-.025em b}\kern-.08em
    T\kern-.1667em\lower.7ex\hbox{E}\kern-.125emX}}
\begin{document}

\title{ Stochastic Gradient Descent  for Constrained Optimization \\
based on Adaptive Relaxed Barrier Functions
% \thanks{Identify applicable funding agency here. If none, delete this.}
}

\author{Naum Dimitrieski, Jing Cao, and Christian Ebenbauer% <-this % stops a space
\thanks{Naum Dimitrieski, Jing Cao, and Christian Ebenbauer are with the Chair of Intelligent Control Systems, RWTH Aachen University, 52062 Aachen, Germany
        {\tt\small (e-mail: \{naum.dimitrieski, christian.ebenbauer\}@ic.rwth-aachen.de, jing.cao@rwth-aachen.de)}
        }%
}

\maketitle

\begin{abstract}

This paper presents a novel stochastic gradient descent algorithm for constrained optimization.
The proposed algorithm randomly samples constraints and components of the finite sum objective function and relies on a relaxed logarithmic barrier function that is appropriately adapted in each optimization iteration.
For a strongly convex objective function and affine inequality constraints, step-size rules and barrier adaptation rules are established
that guarantee asymptotic convergence with probability one. 
The theoretical results in the paper are complemented by numerical studies which highlight potential advantages of the proposed algorithm 
for optimization problems with a large number of constraints.
\end{abstract}
% \begin{IEEEkeywords}
% component, formatting, style, styling, insert
% \end{IEEEkeywords}

\section{Introduction}
Today, stochastic gradient descent (SGD) and its variants are one of the most widely used algorithms in science and engineering. 
SGD-like algorithms are used, for example, in the area of deep learning for training neural networks or in parallel and networked optimization problems. Surprisingly, despite the high relevance of these algorithms, relatively less research about  constrained optimization with SGD has been reported \cite{constrained_SGD_2017good_practice}.
Constraints within the general area of stochastic approximations have been studied, for example, in \cite{kushner2012stochastic}.
SGD-based algorithms based on projection methods have been studied in  \cite{cotter2016light_touch_projection, stoch_recursive_descent_ascent2020objective_subsampling_and_projection, penalization_and_projection_SGD2017admm, constrained_SGD_2017good_practice}, interior point methods in \cite{curtis2023stochastic_interior_point, curtis2024interior_point, he2024stochastic_Interior_Point}, penalty methods in\cite{wang2017stoch_penalty, solis2019continuous_noise_GD_penalty}, and sequential quadratic programming methods in \cite{berahas_curtis2024stochastic, curtis2025almost,qiu2023sequential_quadratic}. All the aforementioned algorithms have in common that their 
iteration steps are based on partial (stochastic) information on the objective function but requiring full (non-stochastic) information on the constraints. 
Only recently, for instance in \cite{xiao2019penalized, gu2024stochastic_projection_Banach_space_barrier_subsampling, penalized_SGD2022new, kervadec2022constrained}, algorithms based on partial information of the constraints have been researched. For example, in \cite{xiao2019penalized, gu2024stochastic_projection_Banach_space_barrier_subsampling}  the authors use partial information of equality constraints to perform an update and projection step. Moreover, in \cite{penalized_SGD2022new} the authors use an adapted softplus penalty function based on partial constraint information at each iteration step. However, the iteration steps of both algorithms rely on projection steps or on solving intermediate optimization problems which can be computationally expensive. Finally, in \cite{kervadec2022constrained}, the authors use an adapted relaxed logarithmic barrier function and combine it in simulation with SGD for optimization of constrained convolutional neural networks. However, no theoretical analysis was presented in \cite{kervadec2022constrained}.

The main contribution of this paper is a novel SGD algorithm for constrained optimization. The algorithm allows for both the objective function and constraints to be only partially known in each iteration and requires no intermediate projection step due to the use of a \textcolor{black}{globally defined} relaxed barrier function. 
In more detail, the contributions of the paper are as follows.  We utilize a relaxed logarithmic barrier function that is appropriately adapted in each iteration \textcolor{black}{and uses} only partial information of the objective function and constraints. For the case of a strongly convex objective function and affine inequality constraints, we derive step size rules and barrier adaption rules that guarantee global convergence with probability one.
Further, we present simulation results which show potential benefits of the proposed algorithm for problems with a large number of constraints. \textcolor{black}{Optimization problems with a large number of constraints are, for example, relevant in deep learning and parallel computing, e.g., \cite{marquez2017DNNhard_constraints,kervadec2022constrained}, \cite{recht2011hogwild}. }

The paper is organized as follows. In \Cref{section_problem_statement} we present the formal problem statement and main results of this paper. We then continue with several numerical results presented in \Cref{section_simulations}, and lastly in \Cref{section_conclusions} we state the main conclusions of this paper and outlook on future works.

\textit{Notation:} The notation is rather standard. $\R^d$ denotes the $d$-dimensional Euclidean space of real numbers and $\R^{+}$ the set of positive real numbers. Moreover, $\E\left[\cdot\right]$ denotes the expectation operator, $\mathcal{C}^n$ denotes the set of $n$-times continuously differentiable functions and a.s. stands for almost surely.
We write a.s. over equations and inequalities to indicate that these relations hold almost surely (with respect to the underlying
probability space).
Lastly, $\overline{a:b}$ and $\mathcal{U}\{\overline{a:b}\}$ for integers~${a,b}$,~${a \leq b}$, denote the set of all integers in the interval $[a,b]$ and a uniform probability distribution (with equal probability) over the integers in $[a,b]$, respectively.

\section{Problem Statement and Main Results}
\label{section_problem_statement}
\begin{figure}[b]
    \centering
    \includegraphics[width=7cm]{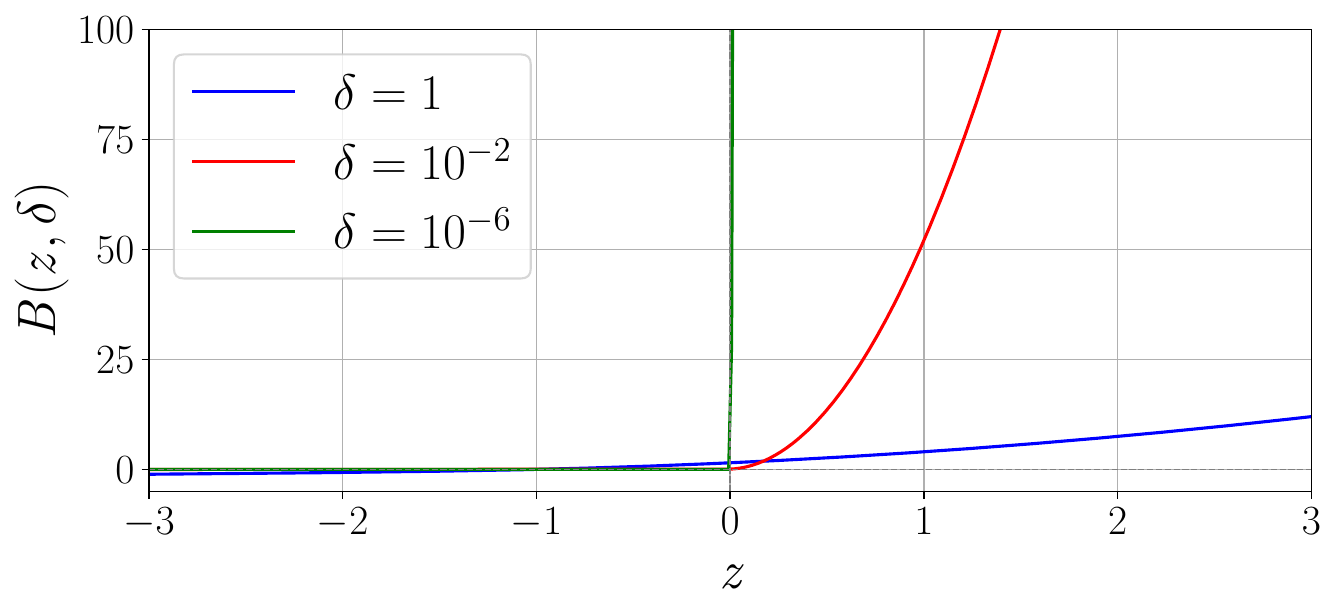}
    \vspace{-0.3cm}
    \caption{Relaxed barrier function \eqref{barrier_fcn} with different values of $\delta > 0$.}
    \label{fig:barrier_fcn_k}
\end{figure}

In this paper, we consider the problem of finding a minimizer of the constrained optimization problem
% \begin{alignat}{2}
%     \label{problem:1}
%     \underset{~}{\mathrm{min}} \hspace{0.23cm}&  f(x) \hspace{1cm} &&\\
%     \mathrm{s.t.} \text{   }& g_j(x) \leq 0,  && j \in \overline{1:m} \notag 
% \end{alignat}
\begin{align}
    \label{problem:1}
    {\mathrm{min}} \hspace{0.1cm}&  f(x) \hspace{0.4cm} \mathrm{s.t.} \hspace{0.4cm} g_j(x) \leq 0, \hspace{0.4cm}  j \in \overline{1:m},
\end{align}
with the finite sum objective function~$${f(x) := \frac{1}{n}\sum_{i=1}^{n}f_i(x)}$$
and with \textcolor{black}{$\mathcal{C}^1$} components~${f_i: \R^d \to \R}$ and \textcolor{black}{$\mathcal{C}^1$} constraints~${g_j:\R^d \to \R}$.
In particular, for our convergence analysis we make the following assumptions.
\textcolor{black}{
\begin{assump}
    \label{assump:1}
    Consider problem \eqref{problem:1}. We assume there exist constants $0 < \mu \leq L$, such that for any $x,y \in \R^d$, and for any $f_i \in \mathcal{C}^1 $, $i \in \overline{1:n}$,
    \begin{align*}
        f_i(y) -  f_i(x) &\leq \nft_i(x)(y-x) + \frac{L}{2}\lVert y-x\rVert_2^2, \\
        f_i(y) -  f_i(x) &   \geq   \nft_i(x)(y-x), \\
        f(y) -  f(x) & \geq \nft(x)(y-x) + \frac{\mu}{2}\lVert y-x\rVert_2^2.
    \end{align*}
\end{assump}
}
%\noindent From \Cref{assump:1} it follows that the  $\underset{x \in \R^d}{\mathrm{min}}f(x) > -\infty$, as $f(x)$ is lower bounded by a quadratic function.
\begin{assump}
    \label{assump:2}
    Consider problem \eqref{problem:1}. We assume there exists a non-empty interior of the feasible set 
    \begin{align*}
        \mathbb{X}:= \{x \in \R^d : g_j(x) \leq 0, \hspace{0.1cm} j \in \overline{1:m}\}.
    \end{align*}
    Moreover, we assume that each $g_j:\R^d \to \R$, $j \in \overline{1:m}$ is affine, i.e., $$g_j(x) := a_j^{\top}x + b_j,$$
    with $a_j \in \R^d$, $b_j \in \R$.
\end{assump}
We denote by $x_C^{*} \in \mathbb{X}$ the minimizer of \eqref{problem:1} under the above assumptions.
\textcolor{black}{Assumption 1 imposes convexity and $L$-smoothness on all $f_i(\cdot)$, $i \in \overline{1:n}$, and strong convexity on $f(\cdot)$. Assumption \ref{assump:1} and \ref{assump:2} are rather restrictive, but not uncommon in 
convergence proofs of gradient-based optimization algorithms \cite{boyd2004convex}. }
The algorithm presented in this paper is based on the idea of barrier methods, i.e., we aim to find a minimizer of \eqref{problem:1} by utilizing a sequence of unconstrained optimization problems, defined for $k \geq 0$ as
\begin{align}
    \label{problem:2}
    \underset{x \in \R^d}{\mathrm{min}} \hspace{0.25cm} f(x) + \frac{1}{m}\sum\limits_{j=1}^{m}B(g_j(x), \delta_k),
\end{align}
where $B: \R \times \R^{+} \to \R$ represents a barrier function adapted by the barrier parameter ${\delta_k \geq \db > 0}$
and $\db$ defines a positive lower bound of the sequence $\{\delta_k\}_{k \ge 0}$. Furthermore, we define
% Furthermore, we define
$${x^{*}(\delta_k) := \underset{x \in \R^d}{\mathrm{argmin}}~f(x) + \frac{1}{m}\sum\limits_{j=1}^{m}B(g_j(x), \delta_k)}.$$ 
Finally, we define a relaxed logarithmic barrier function as used in this paper:
\textcolor{black}{\begin{dfn}
    \label{def:1}
    The function $B: \R \times \R^{+} \to \R$ defined as
    \begin{align}
    \label{barrier_fcn}
        B(z, \delta):= \begin{cases}
            -\delta\log(-z), &  z < - \delta \\
            \frac{1}{2}\left(\frac{(z+2\delta)^2}{\delta} - \delta\right) - \delta\log(\delta), & z \geq -\delta,
        \end{cases}
    \end{align}
    is called a relaxed logarithmic barrier function.
\end{dfn}}
% Finally, for any $\delta > 0$, we define a relaxed logarithm barrier function of the form
% \begin{align}
%     \label{barrier_fcn}
%     B(z, \delta):= \begin{cases}
%         -\delta\log(-z), &  z < - \delta \\
%         \frac{1}{2}\left(\frac{(z+2\delta)^2}{\delta} - \delta\right) - \delta\log(\delta), & z \geq -\delta,
%     \end{cases}
% \end{align}
% which is a convex $\mathcal{C}^2$ map. This barrier function  is defined for all $z \in \R$ and has, for $\delta>0$, a global upper bound on its second derivative. 
\textcolor{black}{This barrier function $B$ is a slightly modified version from the literature (see e.g., \cite{tal1992modified_barrier, hauser2006barrier, feller_ebenbauer2017mpc} for more details on relaxed barrier functions) and differs from the one proposed in \cite{kervadec2022constrained}.}
It is easy to see that $B$ is globally defined, convex in the first argument, and a $\mathcal{C}^2$ map with a $\delta$-dependent global upper bound on its second derivative \cite{feller_ebenbauer2017mpc}. Further,~${\lim_{\db \to 0}B(\cdot, \db)}$ approaches the ideal (limit) barrier function. A graphical illustration of the relaxed barrier function is shown in \Cref{fig:barrier_fcn_k}. We observe that for ${\delta = 1}$ the barrier value increases gradually, while for~${\delta = 0.01}$ the barrier value increases significantly faster. Finally, for~${\delta = 10^{-6}}$ we observe that the shape of the barrier is quite close to that of the ideal barrier, i.e., the barrier 
which has two values, $0$ and $\infty$ for~$z\leq0$ and~$z>0$, respectively.

 The proposed relaxed barrier function offers the advantage over the (non-relaxed) logarithmic barrier function that no projection step or other fallback strategies are needed in order to deal with infeasible iterates, which is specifically challenging for SGD-like algorithms. 
 Furthermore, in contrast to multiplier methods (primal-dual methods), the memory complexity is constant and does not grow with the number of constraints.
 On the other hand, 
 % of course, 
 not all iterates are guaranteed to be feasible when utilizing the proposed relaxation.
 
The sequence $\{x^{*}(\delta_k)\}_{k \geq 0}$ is commonly called the central path 
and intuitively it is plausible that~${\lim_{\db \to 0}x^{*}(\db) = x^{*}_C}$ holds if the interior of the feasible set of \eqref{problem:1} is non-empty (Assumption \ref{assump:2}). We, however, do not give a formal proof in this paper and refer to similar results in the literature \cite[Chapter 11]{boyd2004convex}. 

To find $x_C^{*}$ (respectively $x^{*}(\delta_\infty)$ for $\delta_\infty>0$ close to zero), we propose a stochastic gradient descent algorithm
which simultaneously draws random samples from the objective function and the constraints as follows:
\begin{align}
    \label{SGD}
    x_{k+1} = x_k - \gamma_k \left(\nabla f_{i_k}(x_k) + \nabla B(g_{j_k}(x_k), \delta_k) \right),
\end{align}
where $x_0 \in \R^d$, $i_k \in \overline{1:n}$, $ j_k \in \overline{1:m}$ are i.i.d. random variables,
$\gamma_k$ is the step size and $\delta_k$ is the barrier parameter at iteration $k$. 
%Despite our rather natural approach in form of \eqref{SGD}, we are 
Despite this rather natural Ansatz,
to the best of our knowledge, we are not aware
of any convergence results in the literature for a relaxed barrier-based SGD algorithm as proposed in \eqref{SGD}, \eqref{barrier_fcn}.
For our convergence analysis, we require the following assumption (motivated by the \textcolor{black}{standard} assumption in \cite{robbins1951stochastic}).
\begin{assump}
    \label{assump:3}
    Consider algorithm \eqref{SGD}. We assume that at every $k \geq 0$ it holds that \textcolor{black}{$i_k \distas{} \mathcal{U}\{\overline{1:n}\}$, $  j_k \distas{} \mathcal{U}\{\overline{1:m}\}$.}
\end{assump}
 The uniform distribution has been chosen to ensure that for any $x \in \R^d$ and any~${\delta > 0}$ the standard SGD assumptions hold:
\begin{align*}
        \E\left[\nf_{i_k}(x) \right] &= \frac{1}{n}\sum\limits_{i=1}^{n}\nabla f_i(x) =\nf(x), \\ 
        \E\left[\nb(g_{j_k}(x), \delta) \right] &= \frac{1}{m}\sum\limits_{j=1}^{m}\nb(g_j(x), \delta).
    \end{align*}
\textcolor{black}{In principle, one could choose other distributions, particularly for the barrier indices $j_k$, but this is
beyond the scope of this work.}

We now present our main technical result of this paper which establishes convergence of the iterates generated by \eqref{SGD} to~$x^{*}(\db)$ (respectively arbitrarily close to the minimizer of \eqref{problem:1} for~${\db >0}$ sufficiently small)
with properly tuned step size and barrier adaptation rules.
\begin{thm}
    \label{theorem:convergence}
    Consider problem \eqref{problem:2} and assume that \textcolor{black}{Assumptions} \ref{assump:1} - \ref{assump:3} hold. 
    Let $\{\gamma_k\}_{k \geq 0}$, $\{\ve_k\}_{k \geq 0}$ and  $\{\delta_k\}_{k \geq 0}$ be positive sequences 
    % (i.e., step size rule and adaption rule) 
    such that 
    \begin{align}
    \label{rules}
        \textcolor{black}{a})\sum\limits_{k=0}^{\infty}\gamma_k &= \infty, & \textcolor{black}{b})\sum\limits_{k=0}^{\infty}\gamma_k^2  &< \infty, & \textcolor{black}{c})\sum\limits_{k=0}^{\infty}\gamma_k\ve_k  &< \infty,
    \end{align} with $\delta_k = \db + \ve_k$ and (small) $\db > 0$. 
    % (close to zero)
        Then for any sequence of random variables $\{X_k\}_{k \geq 0}$ defined through algorithm \eqref{SGD}
        it holds that
        \begin{align*}
            \lim\limits_{k \to \infty}X_k \overset{\mathrm{a.s}}{=} x^{*}(\db).
        \end{align*}
\end{thm}
 Due to \Cref{theorem:convergence}, the algorithm iterates converge with probability one based on  conditions on the step size rule
 $\{\gamma_k\}_{k \geq 0}$ which are rather \textcolor{black}{standard \cite{robbins1951stochastic}}.
 The \textcolor{black}{condition $c)$} in \eqref{rules}, however, is not \textcolor{black}{standard} and imposes a condition on the barrier adaptation rule.
Due to the relaxation of the barrier function, the proposed algorithm is not an interior-point algorithm and iterates can be infeasible.
Moreover, in the proof, we require that the relaxed barrier function\textcolor{black}{, in particular, $B(g_j(\cdot),\delta)$, $j \in \overline{1:m}$,} be upper bounded by a quadratic function\textcolor{black}{, which motivates the choice of the barrier function \eqref{barrier_fcn}, as this condition is satisfied in the case of affine inequality constraints}. However, if in \eqref{SGD} we ensure~${\lVert x_k \rVert_2 < \infty}$ for every~${k \geq 0}$ (typically done in the literature via an additional projection step into a compact set), then 
% the statement in
\Cref{assump:2} can be further relaxed, e.g. by requiring all~$g_j(\cdot)$ to be (at least) $\mathcal{C}^1$ convex functions.

%\section{Proof of Convergence}
%\label{section_proof}
In the proof of \Cref{theorem:convergence}, we use the following lemma:
\textcolor{black}{
\begin{lemma}(see \cite[Lemma 1.2]{curtis2025almost}) 
\label{lemma:curtis:robbins}
Let $(\Omega, \Sigma, \mathbb{P})$ be a probability space and let $\{\mathcal{F}_k\}_{k \geq 1}$ with~${\mathcal{F}_k \subseteq \mathcal{F}_{k+1}}$ for all~${k \in \N}$ be a sequence of sub-$\sigma$-algebras of~$\mathcal{F}$. Let~$\{R_k\}_{k \geq 1}$, $\{P_k\}_{k \geq 1}$ and $\{Q_k\}_{k \geq 1}$ be sequences of non-negative random variables such that, for all~${k \in \N}$, the random variables~${R_k,P_k}$ and~$Q_k$ are~$\mathcal{F}_k$ measurable. If~$\sum_{k=1}^{\infty}Q_k < \infty$, and, for all~${k \in \N}$, one has~${\E\left[R_{k+1} | \mathcal{F}_k\right] \leq R_k - P_k + Q_k}$, then, a.s.~${\sum_{k=1}^{\infty}P_k < \infty}$ and~${\lim_{k \to \infty}R_k}$ exists and is finite.
\end{lemma}}

\begin{proof}[Proof of \Cref{theorem:convergence}]
In the following, let~${(\Omega,\Sigma,P)}$ denote a suitable underlying probability space for the considered stochastic process, and let $\{X_k\}_{k\geq0}$, $\{i_k\}_{k\geq0}$, $\{j_k\}_{k\geq0}$ denote sequences of $\Sigma$-measurable random variables where~$\{X_k\}_{k\geq0}$,~${X_k : \Omega \to \R^d}$, is generated by \eqref{SGD}, with~${X_0 = x_0}$, and moreover, ~${i_k: \Omega \to \overline{1:n}}$,~${j_k: \Omega \to \overline{1:m}}$.
Remember that from \Cref{assump:3} we have ${i_k \distas{} \mathcal{U}\{\overline{1:n}\}}$ and ${j_k \distas{} \mathcal{U}\{\overline{1:m}\}}$ at every $k \geq 0$.
For the sake of simplicity, 
we write\textcolor{black}{~${\tx_k: = X_k - x^{*}(\db)}$}
% \begin{align*}
%     \tx_k: = X_k - x^{*}(\db)
% \end{align*}
and~${\lVert \cdot \rVert := \lVert \cdot \rVert_2}$. We denote the natural filtration of the stochastic process $\{X_k,i_k,j_k\}_{k \ge 0}$ defined through the proposed algorithm by~$\mathcal{F}_k$ at time instance~${k \geq 0}$. In addition, we distinguish between relations which hold surely and ones which hold almost surely. We moreover analyze \eqref{SGD} in the following form
\begin{align}
    \label{SGD_adapted}
    X_{k+1} &= X_k - \gamma_k \nabla \Phi_{i_k,j_k}(X_k) + \gamma_k \nabla C_{j_k}(X_k,\delta_k),
\end{align}
with~${\nabla \Phi_{i_k,j_k}(X_k) := \nabla f_{i_k}(X_k) + \nabla B(g_{j_k}(X_k),\db)}$ and ${\nabla C_{j_k}(X_k,\delta_k) :=\nabla B(g_{j_k}(X_k),\db) - \nabla B(g_{j_k}(X_k),\delta_k)}$,
where we denote~${\nabla C_{j_k}(X_k,\delta_k):= \nabla C(g_{j_k}(X_k),\delta_k)}$
and~${\Phi(X_k) := \sum_{i=1}^{n}\frac{ f_i(X_k)}{n} + \sum_{j=1}^{m}\frac{ B(g_j(X_k),\db)}{m}}$ for the sake of simplicity.

The idea of the proof is, \textcolor{black}{first}, to find an upper bound for~${\Ef\left[\lVert \tx_{k+1}  \rVert^2\right] := \E\left[\lVert \tx_{k+1} \rVert^2 | \mathcal{F}_k\right]}$ for every~${k \geq 0}$, and \textcolor{black}{second, to use} this upper bound to show that in the limit~${k \to \infty}$ the sequence converges to~$0$ almost surely. 

\textcolor{black}{\textbf{Step 1:}} Using \eqref{SGD_adapted} and taking expectation conditioned on~$\mathcal{F}_k$, we get
\begin{align*}
    &\Ef\left[\lVert \tx_{k+1} \rVert^2\right] \overset{\mathrm{a.s}}{=} \lVert \tx_k \rVert^2 + \gamma_k^2\Ef\left[\lVert \nabla \Phi_{i_k,j_k}(X_k) \rVert^2\right] \\ +&\gamma_k^2\Ef\left[\lVert \nabla C_{j_k} (X_k,\delta_k) \rVert^2\right] -2\gamma_k\Ef\left[\nabla \Phi_{i_k,j_k}^{\top}(X_k)\tx_k\right] \\ +& 2\gamma_k\Ef\left[\nabla C_{j_k}^{\top} (X_k,\delta_k)\tx_k\right]  \\ -&2\gamma_k^2\Ef\left[\nabla \Phi_{i_k,j_k}^{\top}(X_k)\nabla C_{j_k} (X_k,\delta_k)\right],
\end{align*}
where we remember that $\Ef\left[\tx_k\right] \overset{\mathrm{a.s}}{=} \tx_k$. 
In the following, we provide upper bounds for the right hand side terms. From Cauchy-Schwarz's inequality and Young's inequality we have~${2u^{\top}v \leq 2\|u\|\|v\| \leq \|u\|^2 + \|v\|^2 }$, i.e.,
\begin{align*}
    &-2\gamma_k^2\Ef\left[\nabla \Phi_{i_k,j_k}^{\top}(X_k)\nabla C_{j_k} (X_k,\delta_k)\right] \\
    & \leq \gamma_k^2\Ef\left[\lVert \nabla \Phi_{i_k,j_k}(X_k)\rVert^2 \right] + \gamma_k^2\Ef\left[\lVert \nabla C_{j_k}(X_k, \delta_k)\rVert^2 \right].
\end{align*}
Next, we use Cauchy-Schwarz's inequality to get
\begin{align*}
    \Ef\left[\nabla C_{j_k}^{\top} (X_k,\delta_k)\tx_k\right] & \leq \Ef\left[\lVert \nabla C_{j_k} (X_k,\delta_k) \rVert \lVert \tx_k \rVert \right].
\end{align*}
Furthermore, we use \Cref{assump:3} to get
\begin{align*}
    \Ef\left[\nabla \Phi_{i_k,j_k}^{\top}(X_k)\tx_k\right] &\overset{\mathrm{a.s}}{=} \nabla \Phi^{\top}(X_k)\tx_k,
\end{align*}
where~$\Phi(\cdot)$ is a \textcolor{black}{ strongly convex $\mathcal{C}^1$ function. This 
is the case because $f(\cdot)$ is strongly convex (\Cref{assump:1}) and all $B(g_j(\cdot), \db)$, $j=\overline{1:m}$, are convex (\Cref{assump:2} and \eqref{barrier_fcn}).} 
Thus the inequalities in \cite[Lemma 2.14]{garrigos2023handbook} hold and we get
\begin{align*}
    -2\gamma_k\nabla \Phi^{\top}(X_k)\tx_k &\leq -2\gamma_k\Phi_0(X_k)- \mu\gamma_k \lVert\tx_k\rVert^2,
\end{align*}
where~$\Phi_0(u):=\Phi(u) -  \Phi (x^{*}(\db)) \geq 0$ for all~$u \in \R^d$, as~$x^{*}(\db)$ is the global minimizer of~$\Phi(\cdot)$.
By plugging all the above estimates in~$\Ef\left[\lVert \tx_{k+1} \rVert^2\right]$, we get
\begin{align}
    \label{x_k+1_norm_ineq:1}
     &\Ef\left[\lVert \tx_{k+1} \rVert^2\right] \overset{\mathrm{a.s}}{\leq} \lVert \tx_k \rVert^2 + 2\gamma_k^2\Ef\left[\lVert \nabla \Phi_{i_k,j_k}(X_k) \rVert^2\right] \notag \\ +& 2\gamma_k^2\Ef\left[\lVert \nabla C_{j_k} (X_k,\delta_k) \rVert^2 \right] - \mu\gamma_k \lVert \tx_k\rVert^2 \notag \\ -&2\gamma_k\Phi_0(X_k) + 2\gamma_k \Ef\left[\lVert\nabla C_{j_k} (X_k,\delta_k)\rVert \lVert\tx_k\rVert\right].
\end{align}
In the following we provide further estimates for the right hand side terms using the results from \Cref{appendix}. 

% Starting out with~${\E[\lVert\nabla C_{j_k} (X_k,\delta_k)\rVert \lVert X_k - x^{*}(\db)\rVert | \mathcal{F}_k]}$, from 
As shown in \Cref{appendix A}, \textcolor{black}{where we use \Cref{assump:2}}, we have
\begin{align}
    \label{norm_1_nabla_c_j_avg}
    &\Ef\left[\lVert \nabla C_{j_k}(X_k,\delta_k) \rVert \lVert \tx_k \rVert \right] \overset{\mathrm{a.s}}{\leq}  \lVert \tx_k \rVert^2 \Big( \bar{\bar{a}} \frac{\ve_k}{\delta_k\db} \notag + \C\ve_k\bar{a}  \notag \\+&\bar{b} \frac{\ve_k}{\delta_k\db}\Big)  +\frac{\C\ve_k}{4}\bar{a} +\frac{1}{4}\bar{b} \frac{\ve_k}{\delta_k\db}
\end{align}
with the non-negative constants $\C$, ${\bar{a}, \bar{\bar{a}}, \hat{a}, \bar{b}}$ and $\bar{\bar{b}}$ as defined in \Cref{appendix A}.
Next, as shown in \Cref{appendix B}, \textcolor{black}{where we use \Cref{assump:2}}, we have
\begin{align}
    \label{norm_2_c_j_avg}
    &\Ef\left[\lVert \nabla C_{j_k}(X_k,\delta_k) \rVert^2\right]  \overset{\mathrm{a.s}}{\leq} \textcolor{black}{3\hat{a}\frac{\ve_k^2\lVert \tx_k \rVert^2}{\delta_k^2\db^2}} + 3\C^2\ve_k^2\Bar{\Bar{a}} +
    % \notag
     % \\+& 
     \textcolor{black}{ \frac{3 \Bar{\Bar{b}}\ve_k^2}{\delta_k^2\db^2}}.
\end{align}
Finally, as shown in \Cref{appendix C}, \textcolor{black}{where we use Assumptions~\ref{assump:1} and \ref{assump:2}}, we have
\begin{align}
    \label{norm_Phi_db_2}
    \Ef\left[\lVert \nabla \Phi_{i_k,j_k}(X_k) \rVert^2\right] &\overset{\mathrm{a.s}}{\leq} 4 \Hat{L} \Phi_0(X_k) + 2 \sigma_{\Phi},
\end{align}
with $\Hat{L}: = L + \delta^{-1}\underset{j \in \overline{1:m}}{\mathrm{max}}\lVert a_j \rVert^2$ and $0 \leq \sigma_{\Phi}< \infty$ as defined in \Cref{appendix C}.

We now plug \eqref{norm_1_nabla_c_j_avg}, \eqref{norm_2_c_j_avg} and \eqref{norm_Phi_db_2} into \eqref{x_k+1_norm_ineq:1}, and by algebraic manipulations we get
\begin{align}
    \label{proof:theorem:boundedness_sequence}
    &\Ef\left[\lVert \tx_{k+1} \rVert^2\right] \overset{\mathrm{a.s.}}{\leq} (1-\gamma_k(\mu - \xi_k))\lVert \tx_k \rVert^2 \notag \\ 
    -&2\gamma_k(1-4\gamma_k\Hat{L})\Phi_0(X_k) 
    + \gamma_k\ve_k r_k + 4\gamma^2_k\sigma_{\Phi},
\end{align}
which represents an upper bound estimate with
\begin{align*}
    \xi_k &:= 2 \frac{\ve_k}{\delta_k\db}\left( 3\frac{\gamma_k\ve_k\hat{a}}{\delta_k\db}+ \Bar{\Bar{a}} + \Bar{b} \right) + 2\C\ve_k\bar{a},\\
    r_k &:= \frac{\C\Bar{a}}{2} + \frac{\Bar{b}}{2\delta_k\db} + 6\gamma_k\C^2\ve_k\Bar{\Bar{a}} + 6\gamma_k\frac{\bar{\bar{b}}\ve_k}{\delta_k^2\db^2}.
\end{align*}

In the following, we show that there exists some~${k_0 \geq 0}$, such that for every~${k \geq k_0}$ it holds that~${\mu - \xi_k > 0}$ and~${\gamma_k(1-4\gamma_k\Hat{L}) \geq 0}$. 

    As $\gamma_k$ is bounded and converging to $0$ due to~${\sum_{k=0}^{\infty}\gamma_k^2 < \infty}$, is follows that there exists some~${k_{\gamma} \geq 0}$, such that for all~$k \geq k_{\gamma}$ it holds that~ ${\gamma_k(1-4\gamma_k\Hat{L}) \geq 0}$. Moreover, as~${\sum_{k=0}^{\infty}\gamma_k = \infty}$ and~${\sum_{k=0}^{\infty}\gamma_k\ve_k < \infty}$, it follows that~${\lim_{k\to\infty}\ve_k = 0}$ must be satisfied, where~${\ve_k \in \R^{+}}$
for all~${k \geq 0}$ per assumption. Combined with~${\delta_k \geq \db > 0}$, it follows that~$\xi_k, r_k \in\R^{+}$ for all~${k \geq 0}$, and ~${\lim_{k\to\infty}\xi_k = 0}$. This implies that there exists some time instance~$k_{\xi} \geq 0$, such that for all~$k \geq k_{\xi}$ we have~$\mu  -\xi_k>0$. Finally, we denote~${k_0: =\mathrm{max}(k_{\gamma}, k_{\xi})}$.

We will now show that the iterates $\|X_k\|$ are bounded a.s. for any finite $k$. Remember that~${\|\tx_0\| < \infty}$ by assumption and that~$\xi_k, r_k \in \R^{+}$ for all~${k \geq 0}$. Then, considering \eqref{proof:theorem:boundedness_sequence} and starting from~${k=0}$, each iteration is well-defined and therefore finitely many iterations remain bounded a.s. Since~${k_0 \in \N}$, it follows that~${\|\tx_{k_0}\| \overset{\mathrm{a.s}}{<} \infty}$.

\textcolor{black}{\textbf{Step 2:}} Next, we show the a.s. convergence of the iterates~$X_k$ (for~${k \geq k_0}$) by further analyzing the inequality
\begin{align}
\label{proof:theorem:convergence:prefinal}
    &\Ef\left[\lVert \tx_{k+1} \rVert^2\right] \overset{\mathrm{a.s}}{\leq} (1-\gamma_k(\mu - \xi_k))\lVert \tx_k \rVert^2 \notag \\ 
    +& \gamma_k\ve_k r_k + 4\gamma^2_k\sigma_{\Phi}
\end{align}
for iteration steps ${k \geq k_0}$. For~${k\geq k_0}$, notice that~${-2\gamma_k(1-4\gamma_k\Hat{L})\Phi_0(X_k) \leq 0}$ as $\Phi_0(\cdot)$ is non-negative 
and $\gamma_k(1-4\gamma_k\Hat{L}) \ge 0$
and is therefore not included in the upper bound estimate \eqref{proof:theorem:convergence:prefinal}. Define first~${\Gamma_N:=\sum_{i=k_0}^{N}\gamma_i(\mu - \xi_i) > 0}$ for any $N > k_0$. By using telescopic cancellation 
from ${k = k_0}$ to ${k = N > k_0}$, 
we get
\begin{align*}
    &\E\left[\lVert \tx_{N+1} \rVert^2 | \mathcal{F}_N\right] \overset{\mathrm{a.s}}{\leq} -\sum\nolimits_{k=k_0}^{N}\gamma_k(\mu - \xi_k))\lVert \tx_k \rVert^2 \\ 
    &+ \sum\nolimits_{k=k_0}^{N}\gamma_k\ve_k r_k + 4\sum\nolimits_{k=k_0}^{N}\gamma^2_k\sigma_{\Phi} + \lVert \tx_{k_0} \rVert^2 .
\end{align*}
Rearranging and dividing by $\Gamma_N$
% $\Gamma_N:=\sum_{i=k_0}^{N}\gamma_i(\mu - \xi_i) > 0$, 
we get
\begin{align*}
    &\sum\nolimits_{k=k_0}^{N}p_{N,k}\lVert \tx_k \rVert^2  \overset{\mathrm{a.s}}{\leq} \Gamma_N^{-1}\sum\nolimits_{k=k_0}^{N}\gamma_k(\ve_k r_k + 4\gamma_k\sigma_{\Phi}) \\
    &+\Gamma_N^{-1}\big( \lVert \tx_{k_0} \rVert^2  - \E\left[\lVert \tx_{N+1} \rVert^2 | \mathcal{F}_N\right]\big)   
    =: \Gamma_N^{-1}\bar{M}_N,
\end{align*}
where \textcolor{black}{$p_{N,k} := \gamma_k(\mu - \xi_k)\Big(\sum_{i=k_0}^N\gamma_i(\mu - \xi_i)\Big)^{-1} > 0$}, with the additional property of $\sum_{k=k_0}^{N} p_{N,k} = 1$.
We firstly analyze $\lim_{N\to\infty}\Gamma_N$. As~${\mu - \xi_k > 0}$ for all~$k \geq k_0$ with~${\xi_k > 0}$ by definition, and due to~${\sum_{k=0}^{\infty} \gamma_k = \infty}$ it follows that~$\lim_{N\to\infty}\Gamma_N = \infty$.

We now analyze~$\lim_{N\to\infty} \Bar{M}_N$ and~$\lim_{N\to\infty} \frac{\Bar{M}_N}{\Gamma_N}$. Firstly, observe that~${\lVert \tx_{k_0} \rVert^2  - \E\left[\lVert \tx_{N+1} \rVert^2 | \mathcal{F}_N \right] \leq \lVert \tx_{k_0} \rVert^2 \overset{\mathrm{a.s}}{<} \infty}$.
Next, due to \eqref{rules} and since~${\sigma_{\Phi} < \infty}$ we have~${\sum_{k=0}^{\infty} 4\gamma_k^2\sigma_{\Phi} < \infty}$. Similarly, due to \eqref{rules} and since for every~${k \geq 0}$ it holds~${r_k < \infty}$ we have~${\sum_{k=0}^{\infty} \gamma_k\ve_kr_k < \infty}$. 
We therefore conclude that there exists some~${0 < T  < \infty}$ such that~${\lim_{N\to\infty} \Bar{M}_N \overset{\mathrm{a.s}}{\leq} T <  \infty}$ 
and thus~${\lim_{N\to\infty} \frac{\Bar{M}_N}{\Gamma_N} \overset{\mathrm{a.s}}{\leq} \lim_{N\to\infty} \frac{T}{\Gamma_N} = 0}$. We are then left with\textcolor{black}{~${0 \leq \lim\limits_{N \to \infty}\sum_{k=k_0}^{N}p_{N,k}\lVert \tx_k \rVert^2  \overset{\mathrm{a.s}}{\leq} 0}$}
% \begin{align}
%     \label{proof:theorem:last_eq}
%     0 \leq \lim\limits_{N \to \infty}\sum_{k=k_0}^{N}p_{N,k}\lVert \tx_k \rVert^2  \overset{\mathrm{a.s}}{\leq} 0,
% \end{align}
which \textcolor{black}{further implies} that~$\lim_{N\to\infty}\sum_{k=k_0}^{N}p_{N,k}\lVert \tx_k \rVert^2 \overset{\mathrm{a.s}}{=} 0$. 

We are now left with showing that the sequence~$\{\lVert \tx_k \rVert^2\}_{k \geq k_0}$ converges to zero almost surely. Therefore, we consider \eqref{proof:theorem:convergence:prefinal} for $k \geq k_0$, and refer to \textcolor{black}{\Cref{lemma:curtis:robbins}}. In this setup, we have the non-negative~$\mathcal{F}_k$-measurable random variables~${R_k:=\lVert \tx_k \rVert^2}$,~${P_k := \gamma_k(\mu-\xi_k)\|\tx_k\|^2}$ and~${Q_k: = \gamma_k\ve_k r_k + 4\gamma^2_k\sigma_{\Phi}}$, and we observe that~${\sum_{k=k_0}^{\infty}Q_k < \infty}$. 
Then the result of \textcolor{black}{\Cref{lemma:curtis:robbins}} directly applies, i.e. $\lVert \tx_k \rVert^2$ converges to a finite limit value a.s. 
From the existence of a finite limit value for~$\{\lVert \tx_k \rVert^2\}_{k \geq k_0}$ and\textcolor{black}{~~$\lim_{N\to\infty}\sum_{k=k_0}^{N}p_{N,k}\lVert \tx_k \rVert^2 \overset{\mathrm{a.s}}{=} 0$}, it follows that ${\lim_{N\to\infty}\lVert \tx_k \rVert^2 \overset{\mathrm{a.s}}{=} 0}$,
thus proving the theorem claim.
\end{proof}

\section{Simulations} 
\label{section_simulations}
\begin{figure}[t]
    \centering
    \vspace{0.2cm}
    \includegraphics[width=7.25cm]{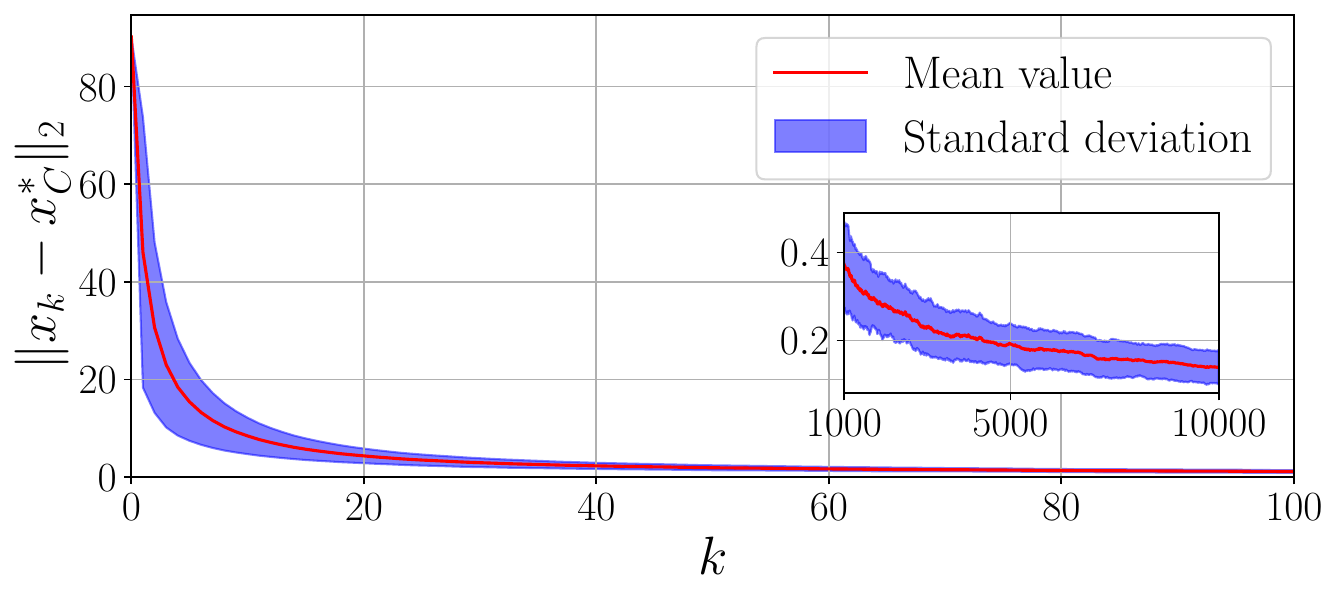}
    \vspace{-0.3cm}
    \caption{Simulation results for an optimization problem with $10^{4}$ constraints. Mean value and estimated standard deviation interval about the mean value are plotted, obtained from a sample set of $1000$ sample trajectories.}
    \label{fig:slow_epsilon_k}
\end{figure}

\begin{figure}[t]
    \centering
    \vspace{0.2cm}
    \includegraphics[width=7.25cm]{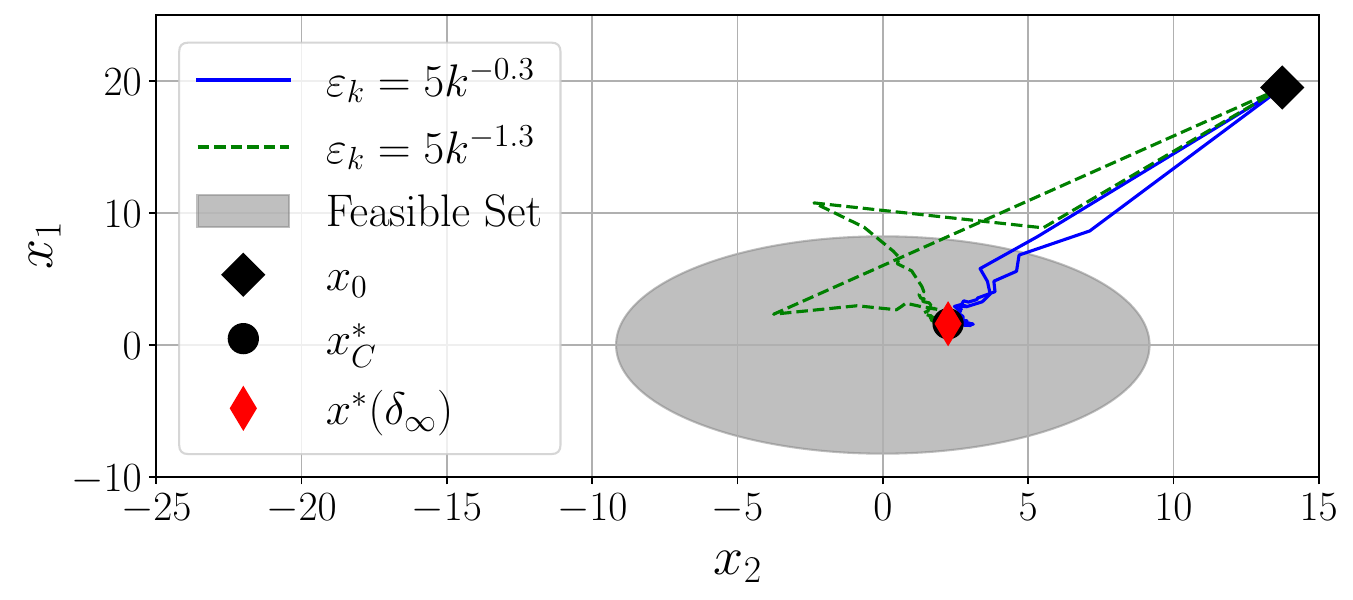}
    \vspace{-0.3cm}
    \caption{Sample trajectories of the proposed SGD algorithm projected onto the $x_1-x_2$ plane. Parameters tuned as in \Cref{fig:slow_epsilon_k} \textcolor{black}{(with $\varepsilon_k = 5k^{-0.3}$ and ${\ve_k = 5k^{-1.3}}$).}}
    \label{fig:compare_epsilon}
\end{figure}

\textcolor{black}{In this section, we apply the proposed adaptive relaxed barrier SGD algorithm, as introduced in  \Cref{section_problem_statement}, to an academic example.}
The presented simulations are carried out in Python (version 3.11.11) using Google Colab on a CPU Intel(R) Xeon(R) @ 2.20GHz, for the following optimization problem. 
The number of decision variables is~${d = 50}$, and the number of inequality constraints varies between~$10$ and~$7\cdot10^{6}$. To setup a feasible set with such a high number of inequality constraints, we consider the ellipsoidal set~${ \mathcal{E}:= \{x \in \R^d: x^{\top}Qx \le 100\}}$, where~$Q=\mathrm{diag}(q_1,...,q_d)$ is a diagonal positive definite matrix
with~${q_i \distas{} \mathcal{U}(1,1.5)}$,~$ {i \in \overline{1:d}}$,
and we randomly sample points from the boundary~${y_j \in \partial\mathcal{E}}$, ${j \in \overline{1:m}}$. 
The affine inequality constraints are then given by~${y_j^{\top}Qx - 100 \leq 0}$ and they define an outer approximation of the ellipsoidal set through the intersection 
of these supporting halfspaces of~$\mathcal{E}$.
This construction ensures that the constrained optimization problem is feasible with arbitrarily many inequality constraints \textcolor{black}{ and that \Cref{assump:2} is satisfied}. We emphasize here that the constraints are sampled once prior to all the simulations. Moreover, the objective function is chosen as
\begin{align*}
    f(x) = \frac{1}{n}\sum\limits_{i=1}^{n}\sum\limits_{j=1}^{d}(\alpha_{i,j} x_j + \log(1 + e^{-\alpha_{i,j}x_j}) + (x_j - \beta)^2)
\end{align*}
with $n=10$, $\alpha_{i,j} \in \R^{+}$ 
and $\beta \in \R$ such that \Cref{assump:1} holds. We have designed the objective function such that~${x_f^{*} = \underset{x \in \R^d}{\mathrm{argmin}}~f(x)} \not \in \mathcal{E}$, i.e., $\lVert x_f^{*} \rVert_2 = 15$. 
At each update step, one of the aforementioned~${n=10}$ functions in~$f(\cdot)$ and one of the $m=10^{4}$ constraints are
sampled to evaluate the gradients of the objective and relaxed barrier functions\textcolor{black}{, thus \Cref{assump:3} is satisfied}.
We first set ${\gamma_k =0.3 k^{-0.8}}$ and~${\db = 10^{-6}}$. For the results presented in \Cref{fig:slow_epsilon_k}
% and \Cref{fig:fast_epsilon_k}
, the barrier \textcolor{black}{parameter is} set as~$\ve_k = 5 k^{-0.3}$
% and~$\ve_k = 5 k^{-1.3}$, respectively
. For this configuration, we obtained~${\lVert x^{*}(\db) - x_C^{*} \rVert_2 = 1.5\cdot 10^{-5}}$. \textcolor{black}{We repeated the same experiment with~$\ve_k = 5 k^{-1.3}$. The results were very similar, but with a slightly increased variance in the transient phase.} 
During the transient phase (from $k=1$ to~approximately~${k=100}$), the variance of the sample trajectories depends on the relative decay of~$\gamma_k$ and~$\ve_k$. 
We \textcolor{black}{have observed} that when $\ve_k$ decreases slower than~$\gamma_k$ 
% (\Cref{fig:slow_epsilon_k}) 
the variance is lower as opposed to when~$\ve_k$ decreases faster than~$\gamma_k$ 
% (\Cref{fig:fast_epsilon_k})
.
However, the asymptotic behavior for both cases seems to be very similar. In addition, in order to illustrate better the convergence of the sample trajectories, in \Cref{fig:compare_epsilon} we present a $2$-D projection of  $4$ sample trajectories, $2$ corresponding to \textcolor{black}{the setup of} \Cref{fig:slow_epsilon_k} (in blue color), and $2$ to the setup \textcolor{black}{of \Cref{fig:slow_epsilon_k} with~$\ve_k = 5 k^{-1.3}$}
% \Cref{fig:fast_epsilon_k} 
(in green color, dashed line). It is worthwhile to mention that choosing a larger $\gamma_k$ and a smaller $\ve_k$ typically leads to a very long transient phase with larger values of $\|x_k-x^{*}(\db)\|_2$.

% \begin{figure}[t]
%     \centering
%     \vspace{0.2cm}
%     \includegraphics[width=8cm]{fig3_fast_epsilon_k.pdf}
%     \caption{Simulation results for an optimization problem with $10^4$ constraints. Parameter $\varepsilon_k$ is chosen to decrease faster than $\gamma_k$. Mean value and estimated standard deviation interval about the mean value are plotted, obtained from a sample set of $1000$ sample trajectories.}
%     \label{fig:fast_epsilon_k}
% \end{figure}

\begin{figure}[t]
    \centering
    \vspace{0.2cm}
    \includegraphics[width=7.5cm]{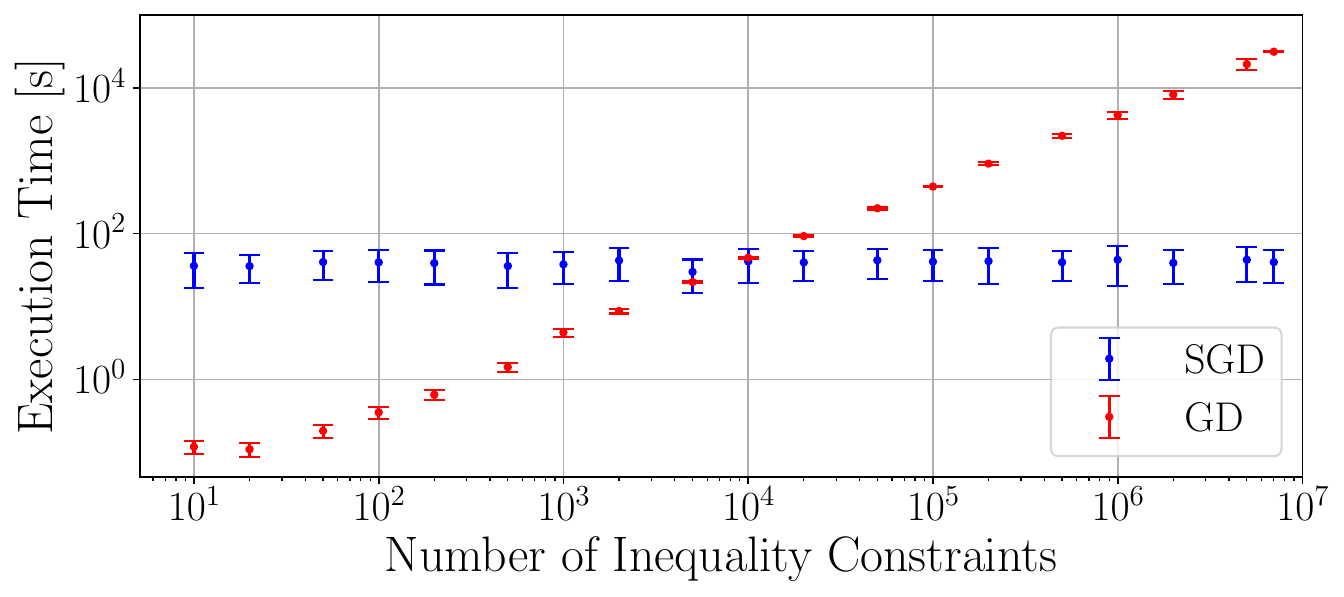}
    \vspace{-0.3cm}
    \caption{Execution time of the proposed SGD algorithm and the deterministic GD algorithm as a function of the number of inequality constraints. 200 sample trajectories were generated for each optimization problem.}
    \label{fig:million_constraints}
\end{figure}

Next, we analyze the execution time of algorithm \eqref{SGD} in terms of the number of constraints, and we compare it to the execution time of the deterministic gradient descent (GD), for which we use the full objective function and barrier constraint information at every update step. The stop criterion for our experiment is the first time instance~$k_\tau$ at which the iterate is close to the constrained optimizer in the sense of~${\lVert x_{k_\tau} - x_C^{*} \rVert_2 \leq 0.01}$.
Notice further that~${\lVert x^{*}(\db) - x_C^{*} \rVert_2}$ is of the order of magnitude~$10^{-5}$ and~\textcolor{black}{${\lVert x_f^{*} - x_C^{*} \rVert_2}$ is of the order of magnitude~$10^{-1}$ and $10^{0}$, however, very similar results have been obtained even for larger values of the latter.}
The sequences are tuned as follows. For algorithm \eqref{SGD} we choose ${\gamma_k= 0.3 k^{-0.8}}$, while for the GD algorithm we choose ${\gamma_k = \gamma = 10^{-2}}$, since for strongly convex problems a constant step size is sufficient to guarantee convergence, and typically leads to a better convergence rate. Further, for both cases we set ${\db = 10^{-6}}$ and~ ${\varepsilon_k = 5 k^{-0.3}}$.
For this setup we present the experiment results in \Cref{fig:million_constraints}. We highlight that the execution time of algorithm \eqref{SGD} remains most often in the range of $40$ to $80$ seconds and is independent of the number of inequality constraints, while the one of the 
% (deterministic) 
GD algorithm grows linearly.
We note further that the break-even point for the execution time is at approximately~${m=10^{4}}$ constraints, and we would like to highlight that at~${m=7\cdot10^{6}}$ the execution time of the GD-based simulations are approximately~$250$ to~$500$ times higher on average than those when using SGD, i.e. algorithm \eqref{SGD}.

\textcolor{black}{In addition, we have used the softplus penalty function from \cite{penalized_SGD2022new} in \eqref{SGD} instead of \eqref{barrier_fcn}, as well as in the aforementioned deterministic GD algorithm. Note that in \cite{penalized_SGD2022new}, the authors do not consider random sampling of constraints at every iteration. The obtained results are very similar to the ones in  \Cref{fig:million_constraints}. Moreover, we implemented a projected gradient descent (PGD) algorithm, where the projection is carried out by \cite[Algorithm 16.3]{nocedal1999numerical}. Algorithm~\eqref{SGD} has an execution time, on average, at least an order of magnitude smaller than the PGD algorithm for~${m=5\cdot10^4}$ constraints, however, the PGD algorithm performed relatively better for up to~${m = 1000}$ constraints, i.e., $0.1$ to $1$ seconds execution time.}

\section{Conclusions and Outlook}
\label{section_conclusions}

In this paper, we have presented a novel barrier-based SGD algorithm. Through a proper barrier relaxation and adaptation rules, the algorithm is guaranteed to converge arbitrarily close to the constrained minimizer of a strongly convex objective function 
% that is 
constrained by affine inequalities while only partial information on the gradients of the objective function and constraints is used in each iteration. As illustrated in  numerical simulations, our approach appears to be more efficient for optimization problems with a large number of inequality constraints than its deterministic counterpart. 

The results in this paper present a first step towards adaptive and relaxed barrier-based SGD algorithms. Many interesting future research questions remain. For example, we have observed that multiple evaluations, in particular $\sqrt{m}$, of the gradients of the barrier functions (constraints) at every iteration step may lead to a decrease in execution time and an increase in convergence rate (as shown for SGD in e.g. \cite{gower2019sgd}).
\textcolor{black}{Moreover, the investigation of post processing procedures 
for the last iterate to achieve feasibility is also an interesting research question.}
Additionally, we believe investigating an iterate-dependent step size and barrier adaptation rule (feedback based, i.e.,~${\gamma(x_k,k),\varepsilon(x_k,k)}$) or non-uniform (adaptive) sampling of the constraints may lead to significant performance increase of the provided algorithm, which could be further complemented with well-known momentum-type acceleration ideas. 
Finally, the convergence result can be improved in terms of less stringent assumptions or explicit convergence rates and an exhaustive benchmarking of the performance of the algorithm is yet to be done.

\bibliographystyle{ieeetr}  
\bibliography{sgd_ref.bib}

\section{Appendices}
\label{appendix}
In the following subsections we prove some inequalities (estimates) that are used in the proof of \Cref{theorem:convergence}. 
Note that throughout this section $X_k$, $\tx_k$ are random variables as in the proof of \Cref{theorem:convergence}, and further, $i_k$ and $j_k$ represent random variables as defined in \Cref{theorem:convergence} that satisfy \Cref{assump:3}.

\subsection{Appendix 1}

\label{appendix 0}
Let $y \in \R^d$, $j \in \overline{1:m}$, and consider $\nabla C_{j}(y,\delta_k)$ as defined for \eqref{SGD_adapted}. Then 
the derivative w.r.t. $y$ (for fixed $\delta_k$) is
\begin{align*}
    \nabla C_{j}(y,\delta_k) &= a_{j} \begin{cases}
       \frac{\delta_k - \db}{g_{j}(y)}, & g_{j}(y) \in \mathcal{I}_1 \\
        \frac{-\db}{g_{j}(y)} - \frac{g_{j}(y) +2\delta_k}{\delta_k}, & g_{j}(y) \in \mathcal{I}_2 \\ 
         g_{j}(y)\frac{\delta_k - \db}{\delta_k\db}, & g_{j}(y) \in \mathcal{I}_3,
    \end{cases}
\end{align*}
where we have $\mathcal{I}_1 := (-\infty, -\delta_k)$, $\mathcal{I}_2 := [-\delta_k, -\db)$ and~${\mathcal{I}_3 := [-\db, \infty)}$. The corresponding second derivative has the form
\begin{align*}
    \nabla^2 C_{j}(y,\delta_k) = a_{j}a_{j}^{\top}\begin{cases}
        -\frac{\delta_k - \db}{g_{j}^2(y)}, & g_{j}(y) \in \mathcal{I}_1 \\
        \frac{\db}{g_{j}^2(y)} - \frac{1}{\delta_k}, & g_{j}(y) \in \mathcal{I}_2 \\ 
         \frac{1}{\db} - \frac{1}{\delta_k}, & g_{j}(y) \in \mathcal{I}_3,
    \end{cases}
\end{align*}
where we observe that the second derivative is negative semi-definite for ~${g_{j}(y) \in \mathcal{I}_1}$, and is positive semi-definite for~${g_{j}(y) \in\mathcal{I}_3}$. This can be easily checked by analyzing the possible values of the second derivative for the respective intervals. Moreover, by analyzing the values of the second derivative for~${g_{j}(y) \in \mathcal{I}_2}$, we get that for~${g_j(y) \in [-\delta_k, -\sqrt{\delta_k\db})}$ it is negative semi-definite and for~${g_{j}(y) \in [-\sqrt{\delta_k\db}, -\db)}$ it is positive semi-definite. It is then easy to verify that~$\|\nabla C_{j}(y,\delta_k)\|$ grows unbounded for~${g_{j}(y) \to \infty}$ and approaches $0$ as~${g_{j}(y) \to -\infty}$. In addition, the reader may verify that for~${g_{j}(y) = -\sqrt{\delta_k\db}}$ we obtain the maximum of the expression~${|\frac{-\db}{g_{j}(y)} - \frac{g_{j}(y) +2\delta_k}{\delta_k}|}$, which equals~${2 - 2\sqrt{\frac{\db}{\delta_k}}}$.
\subsection{Appendix 2}
\label{appendix A}
The goal of this subsection is to find an upper bound of the conditional expectation of~${\lVert \nabla C_{j_k}(X_k, \delta_k) \rVert \lVert \tx_k \rVert}$.
We start out to find an estimate of the form 
\begin{align}
    \label{appendix:2:1}
    \lVert \nabla C_{j_k}(X_k, \delta_k) \rVert & \leq c_0 + \|a_{j_k}\|| g_{j_k}(X_k)| \frac{\delta_k - \db}{\delta_k\db} .
\end{align}
We first consider~${\lVert \nabla C_{j_k}(X_k, \delta_k) \rVert}$ for~${g_{j_k}(X_k) \in \mathcal{I}_3}$ which is~${\|a_{j_k} g_{j_k}(X_k) \frac{\delta_k - \db}{\delta_k\db}\|}$. Applying Cauchy-Schwarz's inequality to this term, we obtain the right hand side in \eqref{appendix:2:1} with $c_0=0$.
In order for this upper bound of~${\lVert \nabla C_{j_k}(X_k, \delta_k) \rVert}$ to be valid for all $g_{j_k}(X_k) \in \R$, we require a sufficiently large constant $c_0$, which we determine in the following.
It is possible to choose
\begin{align*}
    c_0 &: = \mathrm{max}\Big(\underset{g_{j_k}(X_k) \in \mathcal{I}_1}{\mathrm{max}}\| \nabla C_{j_k}(X_k, \delta_k)\|, \\& \hspace{1.43cm}\underset{g_{j_k}(X_k) \in \mathcal{I}_2}{\mathrm{max}}\| \nabla C_{j_k}(X_k, \delta_k)\|\Big) \\
    &=\|a_{j_k}\|\mathrm{max}\Big(\frac{\delta_k - \db}{\delta_k}, 2\Big(1-\sqrt{\frac{\db}{\delta_k}}\Big)\Big),
\end{align*}
where we used that~${\underset{g_{j_k}(X_k) \in \mathcal{I}_1}{\mathrm{sup}}\frac{1}{|g_{j_k}(X_k)|} = \frac{1}{\delta_k}}$, and moreover where~${\underset{g_{j_k}(X_k) \in \mathcal{I}_2}{\mathrm{max}}\|\nabla C_{j_k}(X_k, \delta_k)\| = 2\|a_{j_k}\|\Big(1 - \sqrt{\frac{\db}{\delta_k}}\Big)}$. The latter can be evaluated by analyzing~$\|\nabla C_{j_k}(X_k, \delta_k)\|$ over the interval~${g_{j_k}(X_k) \in \mathcal{I}_2}$, where according to \Cref{appendix 0} ${g_{j_k}(X_k) = -\sqrt{\db\delta_k}}$ is the maximizer 
of~${|\frac{-\db}{g_{j_k}(X_k)} - \frac{g_{j_k}(X_k) +2\delta_k}{\delta_k}|}$ on said interval, which we leave to the reader to verify. Then, using $\delta_k = \db + \ve_k$ we reformulate 
%\begin{align*}
    $1-\sqrt{\frac{\db}{\delta_k}} = \frac{1 - \frac{\db}{\delta_k}}{1+\sqrt{\frac{\db}{\delta_k}}} = \frac{\ve_k}{\delta_k + \sqrt{\db\delta_k}}$
%\end{align*}
to obtain
\begin{align*}
    c_0 &= \ve_k\|a_{j_k}\| \mathrm{max}\Big(\frac{1}{\delta_k}, \frac{2}{\delta_k + \sqrt{\delta_k\db}}\Big)
    \\ &\leq \ve_k\|a_{j_k}\| \underset{k \geq 0}{\mathrm{max}}\Big(\mathrm{max}\Big(\frac{1}{\delta_k}, \frac{2}{\delta_k + \sqrt{\delta_k\db}}\Big)\Big) \\
    &=:\ve_k\|a_{j_k}\|\C.
\end{align*}
Then, by plugging in $\delta_k = \db + \ve_k$, we get
\begin{align*}
    \lVert \nabla C_{j_k}(X_k,\delta_k)\rVert &\leq \C\ve_k\|a_{j_k}\| + \lVert a_{j_k} \rVert | g_{j_k}(X_k) | \frac{\ve_k}{\delta_k\db}.
\end{align*}
Next, expressing~${g_{j_k}(X_k) = a_{j_k}^{\top}\tx_k + b_{j_k} + a_{j_k}^{\top}x^{*}(\db)}$, and thereafter using~${\lVert u + v \rVert \leq \lVert u \rVert + \lVert v \rVert}$ and Cauchy-Schwarz's inequality we obtain that
\begin{align}
    \label{nabla_c_j_rearranged}
    &\lVert \nabla C_{j_k}(X_k,\delta_k)\lVert \leq \lVert a_{j_k}\rVert \frac{\ve_k}{\delta_k\db}  \Big(\C\delta_k\db+\lVert a_{j_k}\rVert \lVert\tx_k \rVert \notag \\ +&|b_{j_k} + a_{j_k}^{\top}x^{*}(\db) |\Big) =:  D(X_k,\delta_k).
\end{align}
Further, by reformulation we have
\begin{align*}
    &\lVert \nabla C_{j_k}(X_k,\delta_k) \rVert \lVert \tx_k \rVert \leq \lVert a_{j_k} \rVert^2 \frac{\ve_k}{\delta_k\db}\lVert \tx_k \rVert^2 \\ +& \C\ve_k \lVert a_{j_k} \rVert \lVert \tx_k \rVert \notag  + \lVert a_{j_k} \rVert | b_{j_k} + a_{j_k}^{\top}x^{*}(\db)| \frac{\ve_k}{\delta_k\db}\lVert \tx_k \rVert.
\end{align*}
We now want to upper bound the right hand side by an affine function in $\lVert \tx_k\rVert^2$ as utilized in the proof of \Cref{theorem:convergence}.  Therefore, we use the inequality between the arithmetic and geometric mean to obtain~${    \frac{1}{4} + \lVert \tx_k\rVert^2 \geq \lVert \tx_k \rVert}$,
which we plug in the above equation to get
\begin{align*}
    % \label{norm_1_nabla_c_j}
    &\lVert \nabla C_{j_k}(X_k,\delta_k) \rVert \lVert \tx_k \rVert \leq \lVert \tx_k \rVert^2 \Big( \lVert a_{j_k} \rVert^2 \frac{\ve_k}{\delta_k\db} \notag \\+& \C\ve_k\lVert a_{j_k} \rVert  \notag
    +\lVert a_{j_k} \rVert| b_{j_k}+ a_{j_k}^{\top}x^{*}(\db) | \notag \frac{\ve_k}{\delta_k\db}\Big) \notag +\frac{\C\ve_k\lVert a_{j_k} \rVert }{4} \notag \\
    +&\frac{1}{4}\lVert a_{j_k} \rVert| b_{j_k} + a_{j_k}^{\top}x^{*}(\db) |  \frac{\ve_k}{\delta_k\db}.
\end{align*}
Next, we will have a look at 
the $\mathcal{F}_k$-conditioned expectation of the above term. As~${j_k \distas{} \mathcal{U}\{\overline{1:m}\}}$, it follows that 
\begin{align*}
    \Ef\left[\lVert \nabla C_{j_k}(X_k,\delta_k) \rVert \lVert \tx_k \rVert\right] \overset{\mathrm{a.s}}{=} \sum_{j=1}^{m}\frac{\lVert \nabla C_j(X_k,\delta_k) \rVert \lVert \tx_k \rVert}{m},
\end{align*}
where we define ${\Bar{a}:=\frac{1}{m}\sum_{j=1}^{m} \lVert a_j \rVert}$, ${\Bar{\Bar{a}}:=\frac{1}{m}\sum_{j=1}^{m} \lVert a_j \rVert^2}$, ${\hat{a}:=\frac{1}{m}\sum_{j=1}^{m} \lVert a_j \rVert^4}$, as well as the terms~${\Bar{b}:=\frac{1}{m}\sum_{j=1}^{m} \lVert a_j \rVert | b_j + a_j^{\top}x^{*}(\db)|}$ and ${\Bar{\Bar{b}}:=\frac{1}{m}\sum_{j=1}^{m} \lVert a_j \rVert^2 | b_j + a_j^{\top}x^{*}(\db) |^2}$.
Thus, we obtain
\begin{align}
    &\Ef\left[\lVert \nabla C_{j_k}(X_k,\delta_k) \rVert \lVert \tx_k \rVert \right] \overset{\mathrm{a.s}}{\leq}  \lVert \tx_k \rVert^2 \Big( \bar{\bar{a}} \frac{\ve_k}{\delta_k\db} \notag + \C\ve_k\bar{a}  \notag \\+&\bar{b} \frac{\ve_k}{\delta_k\db}\Big) \notag +\frac{\C\ve_k}{4}\bar{a} +\frac{1}{4}\bar{b} \frac{\ve_k}{\delta_k\db}. \notag
\end{align}

\subsection{Appendix 3}
\label{appendix B}
Consider the term $\Ef\left[\lVert \nabla C_{j_k}(X_k,\delta_k) \rVert^2\right]$ as given in \Cref{appendix 0}, and  remember that \eqref{nabla_c_j_rearranged} applies. Then we also have
\begin{align*}
    \lVert \nabla C_{j_k}(X_k,\delta_k) \rVert^2 \leq D^2(X_k,\delta_k),
\end{align*}
and by using~${\lVert u + v + y\rVert^2 \leq 3\lVert u \rVert^2 + 3\lVert v \rVert^2 + 3\lVert y \rVert^2}$ we get 
\begin{align*}
    &\lVert \nabla C_{j_k}(X_k,\delta_k) \rVert^2  \leq 3\lVert a_{j_k} \rVert^4\frac{\ve_k^2}{\delta_k^2\db^2}\lVert \tx_k \rVert^2 + 3\C^2\ve_k^2\lVert a_{j_k} \rVert^2
    \notag \\  +&3 \lVert a_{j_k} \rVert^2 | b_{j_k} + a_{j_k}^{\top}x^{*}(\db) |^2 \frac{\ve_k^2}{\delta_k^2\db^2}.
\end{align*}
As shown in \Cref{appendix A} for \eqref{norm_1_nabla_c_j_avg}, 
we use the conditioned expectation on $\mathcal{F}_k$ and \Cref{assump:3} to arrive at the desired estimate
\begin{align}
    &\Ef\left[\lVert \nabla C_{j_k}(X_k,\delta_k) \rVert^2\right]  \overset{\mathrm{a.s}}{\leq} 3\hat{a}\frac{\ve_k^2}{\delta_k^2\db^2}\lVert \tx_k \rVert^2 + 3\C^2\ve_k^2\Bar{\Bar{a}} \notag
     \\+3& \Bar{\Bar{b}} \frac{\ve_k^2}{\delta_k^2\db^2}. \notag
\end{align}

\subsection{Appendix 4}
\label{appendix C}

The goal of this subsection is to verify that the first inequality in \cite[Lemma 2.29]{garrigos2023handbook} 
about convex $L$-smooth functions holds for every $\Phi_{i_k,j_k}(\cdot)$, i.e., for every $i_k,j_k$, and furthermore to derive an upper bound  for~${\Ef\left[\lVert \nabla \Phi_{i_k,j_k}(X_k) \rVert^2\right]}$.

We first require some preliminary analysis. We start by considering the functions $\Phi_{i_k,j_k}(\cdot)$ and $\Phi(\cdot)$ as defined in \eqref{SGD_adapted}. 
We remember that~${\Phi_{i_k,j_k}(\cdot) = f_{i_k}(\cdot) + B(g_{j_k}(\cdot),\db)}$, where by \Cref{assump:1} ${f_{i_k}(\cdot) \in \mathcal{C}^1}$ and is convex, and by design and \Cref{assump:2} ${B(g_{j_k}(\cdot),\db) \in \mathcal{C}^2}$ and is convex in the first argument. Hence, $\Phi_{i_k,j_k}(\cdot) \in \mathcal{C}^1$ and is convex, and this holds for any realization of $i_k,j_k$.
Similarly, we observe that, by definition,~${\Phi(\cdot) = f(\cdot) + \sum_{j=1}^{m}\frac{ B(g_j(\cdot),\db)}{m}}$. From \Cref{assump:1} 
we have that ${f(\cdot) \in \mathcal{C}^1}$ is a strongly convex function. Using the fact that all~${B(g_{j}(\cdot),\db) \in \mathcal{C}^2}$,~${j \in \overline{1:m}}$, are convex, it follows that ${\Phi(\cdot) \in \mathcal{C}^1}$ is a strongly convex function. In addition, we remember that~${\nabla \Phi(x^{*}(\db)) = 0}$.
Moreover, from the definition of the barrier function in \eqref{barrier_fcn} we get
\begin{align}
    \label{nabla_2_B_upper_bound}
    \nabla^2 B(g_{j_k}(x), \db) &= a_{j_k}a_{j_k}^{\top}\begin{cases}
        \frac{\db}{(g_{j_k}(x))^2}, & g_{j_k}(x) < -\db \\
        \frac{1}{\db}, & g_{j_k}(x) \geq -\db
    \end{cases} \notag \\
    & \preceq \db^{-1} \underset{j \in \overline{1:m}}{\mathrm{max}} \lVert a_{j} \rVert^2 I,
\end{align}
and we define $\Hat{L}:= L + \db^{-1} \underset{j \in \overline{1:m}}{\mathrm{max}} \lVert a_{j} \rVert^2$.

Next, we verify that the first inequality in \cite[Lemma 2.29]{garrigos2023handbook} holds for every $\Phi_{i_k,j_k}(\cdot)$, and we start by considering any~${\theta_1,\theta_2, \theta_3 \in \R^d}$. By definition, we have
\begin{align*}
    \Phi_{i_k,j_k}(\theta_2) = f_{i_k}(\theta_2) + B(g_{j_k}(\theta_2),\db).
\end{align*}
Using \Cref{assump:1} on $f_{i_k}(\theta_2)$ we further get
\begin{align*}
    f_{i_k}(\theta_2) & \leq f_{i_k}(\theta_1) + \nabla f_{i_k}^{\top}(\theta_1)(\theta_2 - \theta_1) + \frac{L}{2}\lVert \theta_2 - \theta_1\rVert^2.
\end{align*}
In addition, using the Taylor expansion of~${B(g_{j_k}(\cdot),\db) \in \mathcal{C}^2}$ and the upper bound of~$\nabla^2 B(g_{j_k}(\cdot),\db)$, as defined in \eqref{nabla_2_B_upper_bound}, we obtain
\begin{align*}
    &B(g_{j_k}(\theta_2),\db)  \leq B(g_{j_k}(\theta_1),\db) \\ +& \nabla B^{\top}(g_{j_k}(\theta_1),\db)(\theta_2 - \theta_1) + \frac{\db^{-1}\underset{j \in \overline{1:m}}{\mathrm{max}}\lVert a_j \rVert^2}{2}\lVert \theta_2 - \theta_1\rVert^2.
\end{align*}
We combine the last two inequalities and use the definition of~$\Phi_{i_k,j_k}(\cdot)$. Thereafter, we employ elementary algebraic transformations to obtain
\begin{align}
    \label{appendix:3:1}
    \Phi_{i_k,j_k}(\theta_2) - \Phi_{i_k,j_k}(\theta_1) &\leq \nabla \Phi_{i_k,j_k}^{\top}(\theta_1)(\theta_2 - \theta_1) \notag \\ &+ \frac{\Hat{L}}{2}\lVert \theta_2 - \theta_1\rVert^2.
\end{align}
Similarly, as $\Phi_{i_k,j_k}(\cdot)$ is convex, we have
\begin{align}
    \label{appendix:3:2}
    \Phi_{i_k,j_k}(\theta_3) - \Phi_{i_k,j_k}(\theta_2) &\leq \nabla \Phi_{i_k,j_k}^{\top}(\theta_3)(\theta_3 - \theta_2). 
\end{align}
Summing \eqref{appendix:3:1} and \eqref{appendix:3:2} gives us 
\begin{align}
    \label{appendix:3:3}
   & \Phi_{i_k,j_k}(\theta_3) - \Phi_{i_k,j_k}(\theta_1) \leq \nabla \Phi_{i_k,j_k}^{\top}(\theta_3)(\theta_3 - \theta_2) \notag \\ &+  \nabla \Phi_{i_k,j_k}^{\top}(\theta_1)(\theta_2 - \theta_1) + \frac{\Hat{L}}{2}\lVert \theta_2 - \theta_1\rVert^2.
\end{align}
Observe that by setting~$\theta_3 = x$, ${\theta_2 = z}$ and~$\theta_1 = y$, we obtain the last inequality on page $12$ in the proof of \cite[Lemma 2.29]{garrigos2023handbook}. We leave it to the reader to observe that the follow-up analysis presented in said proof applies to $\Phi_{i_k,j_k}(\cdot)$. We thereby verify that the first inequality in \cite[Lemma 2.29]{garrigos2023handbook} holds for $\Phi_{i_k,j_k}(\cdot)$, which ensures that the following holds
\begin{align}
    \label{appendix:4:1}
     &\lVert \nabla \Phi_{i_k,j_k}(X_k) - \nabla \Phi_{i_k,j_k}(x^{*}(\db)) \rVert^2 \leq 2\Hat{L}\Big( \Phi_{i_k,j_k}(X_k)  \notag \\ -& \Phi_{i_k,j_k}(x^{*}(\db))  - \nabla \Phi_{i_k,j_k}^{\top}(x^{*}(\db))\tx_k\Big).
\end{align}
We note here that the constant~$\Hat{L}$ remains the same for every realization of~$i_k$ and~$j_k$, i.e. for every~$\Phi_{i_k,j_k}(\cdot)$.

We now proceed with deriving the upper bound estimate for~$\Ef\left[\| \nabla \Phi_{i_k,j_k}(X_k) \|^2\right]$. Similar to  \cite[Lemma 4.20]{garrigos2023handbook}, we start out by rewriting
\begin{align*}
    &\lVert \nabla \Phi_{i_k,j_k}(X_k) \rVert^2 = \lVert \nabla \Phi_{i_k,j_k}(X_k) -  \nabla \Phi_{i_k,j_k}(x^{*}(\db)) \\+& \nabla \Phi_{i_k,j_k}(x^{*}(\db)) \rVert^2, 
\end{align*}
and using the inequality ${\lVert u + v\rVert^2 \leq 2\lVert u \rVert^2 + 2\lVert v \rVert^2 }$ to get
\begin{align*}
    &\lVert \nabla \Phi_{i_k,j_k}(X_k) \rVert^2 \leq 2\lVert \nabla \Phi_{i_k,j_k}(X_k) -  \nabla \Phi_{i_k,j_k}(x^{*}(\db)) \rVert^2 \\ +& 2\lVert \nabla \Phi_{i_k,j_k}(x^{*}(\db)) \rVert^2. 
\end{align*}
Now, using \cite[Lemma 2.29]{garrigos2023handbook} in the above expression gives 
\begin{align*}
    &\lVert \nabla \Phi_{i_k,j_k}(X_k) \rVert^2 \leq 4\Hat{L}\Big(\Phi_{i_k,j_k}(X_k) -  \Phi_{i_k,j_k}(x^{*}(\db)) \\ -&\nabla \Phi_{i_k,j_k}^{\top}(x^{*}(\db))(X_k - x^{*}(\db))\Big)  \\ +& 2\lVert \nabla \Phi_{i_k,j_k}(x^{*}(\db)) \rVert^2.
\end{align*}
Taking the conditioned expectation on $\mathcal{F}_k$, we get
\begin{align*}
    &\Ef\left[\lVert \nabla \Phi_{i_k,j_k}(X_k) \rVert^2\right] \overset{\mathrm{a.s}}{\leq} 4 \Hat{L} \Big( \Ef\left[\Phi_{i_k,j_k}(X_k)\right] \\ -& \Ef\left[\Phi_{i_k,j_k}(x^{*}(\db))\right]\Big)
      -4\hat{L}\Ef\left[\nabla \Phi_{i_k,j_k}^{\top}(x^{*}(\db))\tx_k\right]
    \\  +& 2 \Ef\left[\lVert \nabla \Phi_{i_k,j_k}(x^{*}(\db)) \rVert^2\right],
\end{align*}
where, by \Cref{assump:3} and the definition of $\Phi(\cdot)$, we have 
\begin{align*}
    \Ef\left[\nabla \Phi_{i_k,j_k}^{\top}(x^{*}(\db))\tx_k\right] &\overset{\mathrm{a.s}}{=} \nabla \Phi^{\top}(x^{*}(\db))\tx_k= 0 \\
    \Ef\left[\Phi_{i_k,j_k}(X_k)\right] &\overset{\mathrm{a.s}}{=} \Phi(X_k) \\
    \Ef\left[\Phi_{i_k,j_k}(x^{*}(\db))\right] &= \Phi(x^{*}(\db)) = \underset{x \in \R^d}{\mathrm{min}} \Phi(x).
\end{align*}
Finally, we define ${\sigma_{\Phi}:= \Ef\left[\lVert \nabla \Phi_{i_k,j_k}(x^{*}(\db)) \rVert^2\right]}$. Now we show that~$\sigma_{\Phi}$ is well-defined and bounded, and using \Cref{assump:3}, we get
\begin{align}
    \label{appending:4:final}
    \sigma_{\Phi} = \frac{1}{nm}\sum_{i=1}^{n}\sum_{j=1}^{m}\lVert \nabla \Phi_{i,j}(x^{*}(\db)) \rVert^2.
\end{align}
Remember that for every realization of $i_k,j_k$ the function~$\Phi_{i_k,j_k}(\cdot) \in \mathcal{C}^1$, thus the gradient is globally well-defined. Moreover, since $\Phi(\cdot)$ is strongly convex,~${x^{*}(\db) \in \R^d}$ exists. It is then evident that each term in the sum of \eqref{appending:4:final} is a non-negative real number, thereby proving that $\sigma_{\Phi} < \infty$.
Thus, we get
\begin{align}
    \Ef\left[\lVert \nabla \Phi_{i_k,j_k}(X_k) \rVert^2\right] &\overset{\mathrm{a.s}}{\leq} 4 \Hat{L} ( \Phi(X_k) - \Phi(x^{*}(\db))) \notag \\ &+ 2 \sigma_{\Phi}. \notag
\end{align}

\end{document}